\documentclass[11pt,twoside,leqno,letterpaper]{article}
\usepackage{amsmath,amssymb}
\makeatletter
\newcommand{\provsection}{\@startsection%
{section}
{1}
{\parindent}
{-\baselineskip}
{-0.5\baselineskip}
{\normalfont\Large\bfseries}
}%
\makeatother%
\renewcommand{\section}[1]{\provsection{#1}}
\makeatletter
\newcommand{\provsubsection}{\@startsection%
{subsection}
{2}
{\parindent}
{-\baselineskip}
{-0.5\baselineskip}
{\normalfont\normalsize\bfseries}
}%
\makeatother%
\renewcommand{\subsection}[1]{\provsubsection{#1.}}
\makeatletter
\newcommand{\provsubsubsection}{\@startsection%
{subsubsection}
{3}
{\parindent}
{-0.2\baselineskip}
{-0.5\baselineskip}
{\normalfont\normalsize\bfseries}
}%
\makeatother%
\renewcommand{\subsubsection}[1]{\provsubsubsection{#1}}
 
\renewcommand{\[}{\begin{eqnarray*}}
\renewcommand{\]}{\end{eqnarray*}}
\newcommand{\la}{\begin{eqnarray}}
\newcommand{\al}{\end{eqnarray}}

\newcommand{\conv}{\ast}

\DeclareMathOperator*{\bigconv}{\mbox{\LARGE$\conv$}}
\newcommand{\Halmos}
{\mbox{}\hfill$\square$}

\renewcommand{\d}{{\,\text{\rm d}}}

\newcommand{\set}[1]{{\left\{#1\right\}}}

\renewcommand{\epsilon}{\varepsilon}\renewcommand{\phi}{\varphi}      %
\renewcommand{\rho}{\varrho}        \renewcommand{\theta}{\vartheta}  %
\newcommand{\N}{{\mathbb N}}
\newcommand{\R}{{\mathbb R}}
\newcommand{\Z}{{\mathbb Z}}

\renewcommand{\P}{{\mathbb P}}                                               %
\newcommand{\given}{\,\pmb{\big\vert}\,}                                     %
\newcommand{\Ber}{\mathrm{B}}     
\newcommand{\BerC}{\mathrm{B}}    
\newcommand{\STPC}{\mathrm{Q}}    
\newcommand{\stge}{\ge_{\text{st}}}     

\newcommand{\pollard}{ \sim\kern -0.7ex>}

\newcommand{\AnnMS}{{\it Ann.\ Math.\ Statist. }}               %
\newcommand{\BAMS}{{\it Bull.\ Amer.\ Math.\ Soc. }} 
\newcommand{\SPL}{{\it Statistics \&  Probab.\ Letters }}       %
\newcommand{\TPA}{{\it Theory Probab.\ Appl. }}                 %

\title{Lower bounds for tails  of sums \\
of independent symmetric random variables }

\author{Lutz Mattner,\\
         Universit\"at zu L\"ubeck   }

\begin{document}
\maketitle
\begin{abstract}
The approach  of Kleitman (1970) and Kanter (1976)
to multivariate concentration function inequalities is
generalized in order to obtain for deviation probabilities of sums
of independent symmetric random variables a lower bound depending
only on deviation probabilities of the terms of the sum.
This bound is optimal up to discretization effects, improves on a result of
Nagaev (2001), and complements the comparison theorems of 
Birnbaum (1948) and Pruss (1997). 
Birnbaum's theorem for unimodal random variables
is extended to the lattice case.
\end{abstract}

\footnotetext{{\em 2000 Mathematics Subject Classification:} 
 60E15,   60G50. 
 }
\footnotetext{{\em Keywords and phrases:} 
 Bernoulli convolution,
 concentration function,
 deviation probabilities,
 Poisson binomial distribution,
 symmetric three point convolution,
 unimodality.}


\section{Introduction} \mbox{} \\

For deviation probabilities $\P (|S| >t)$ of sums
$S = \sum_{i=1}^n X_i$ of independent, real-valued,
and symmetrically distributed  random variables $X_i$,
Nagaev (2001, Theorem 1, in different notation)
obtained the lower bound
\la \label{Nagaev Thm 1}
  \P (|S| >t) &\ge & \sum_{k>t/h} 2^{-k}\mathrm{B}_p(\set{k}) 
 \qquad \qquad (t \in [0,nh[)
\al
where $h\in\,]0,\infty[$ is a free parameter and 
\la   \label{Def BerC}
 \BerC_p &:=& \bigconv_{i=1}^n \Ber_{p_i}
\al
is the convolution of the Bernoulli distributions 
$\Ber_{p_i}=(1-p_i)\delta_{0} +p_i\delta_1$ with 
success probabilities $p_i := \P(|X_i| \ge h)$.
Nagaev also provided analytically more tractable  lower bounds 
for the right hand side of \eqref{Nagaev Thm 1}
and showed that the resulting inequalities for $ \P (|S| >t)$
effectively complement other bounds depending on  second and  third 
absolute moments of the random variables $X_i$.

The main purpose of the present note is to provide 
as Theorem \ref{Theorem (Kanter's (1976) Lemma 4.2 generalized)} 
below a generalization of Kanter's (1976) concentration function inequality
for sums of independent and symmetric random vectors, which yields as 
Corollary \ref{Corollary improving Nagaev} below  
in particular the following improvement of \eqref{Nagaev Thm 1},
under the same assumptions as above:
\la  \label{Improvement of Nagaev Thm 1} 
 \P (|S| >t) &\ge & 
  \sum_{k>t/h} \Big(1- 2^{-k} F_k(\lfloor \frac th +1\rfloor) \Big)
  \mathrm{B}_p(\set{k})     \qquad \qquad (t \in [0,nh[)
\al 
Here and below, we use the standard notations
$\lfloor x\rfloor :=  \max\set{k\in\Z : k\le x   }$ 
and $\lceil x \rceil := -\lfloor -x\rfloor$, and write
\la \label{Def Fnk, sums of largest binomial coefficients}
 F_n(m) &:=& \max_{r\in\Z} \sum_{i=r}^{r+m-1} {n\choose i} 
 \qquad  \qquad (n,m\in \N_0)
\al
for the sum of the $m$ largest binomial coefficients of 
order $n$. 
For $t \in [(n-1)h ,nh [$, the inequalities in 
\eqref{Nagaev Thm 1} and \eqref{Improvement of Nagaev Thm 1} 
are identical, while for $t \in [0, (n-1)h[$
and $\BerC_p(\set{n-1})>0$, inequality  
\eqref{Improvement of Nagaev Thm 1} is strictly 
sharper than \eqref{Nagaev Thm 1}. 
Moreover, as follows from the proof of Corollary 
\ref{Corollary improving Nagaev},
inequality \eqref{Improvement of Nagaev Thm 1} 
is optimal up to discretization effects, in the sense that,
subject to the stated assumptions, the right hand side of  
\eqref{Improvement of Nagaev Thm 1} is the greatest lower
bound for $\P(|S| >t) +\frac 12 \P(|S|=t)$
for every $t=mh$ with $m\in\set{1,\ldots,m}$.

The rest of this note is structured as follows. 
Section \ref{A generalized Kanter inequality} develops the 
Kleitman-Kanter approach to multivariate concentration function
inequalities. A specialization to the one-dimensional 
case, namely Corollary \ref{Corollary to Kanter Lemma 4.2, Version 2},
immediately yields the above-mentioned Corollary 
\ref{Corollary improving Nagaev} improving  Nagaev's result. 
Section \ref{Comparison theorems} reformulates 
Corollary \ref{Corollary to Kanter Lemma 4.2, Version 2}
as a comparison theorem, stated together with  related
results of Pruss (1997) and Birnbaum (1948). 
The latter is generalized to the lattice case. 
Historical remarks are collected in Section 
\ref{Historical notes}.

\section{A generalized Kanter inequality}
\label{A generalized Kanter inequality}
 \mbox{} \\

Let $\|\cdot\|$ be a seminorm on an $\R$-vector-space $E$
and let $|\cdot|$ denote the usual absolute value on $\R$.
We write $\N:= \set{1,2,3,\ldots}$ and $\N_0 := \set{0}\cup\N$.

\subsection{Lemma} \label{Kleitman-Lemma}  
{\sl Let $a\in E$, $m\in\N$ and
$C_1,\ldots,C_m\subset E$ with
\la          \label{Kleitman-Lemma assumption}
 x,y \in C_j &\Rightarrow& \|x-y\| < \| a \|
\al 
for each $j\in\set{1,\ldots,m}$. Then for some $r\in\set{1,\ldots,m}$
the translate  $C_r-a$ is disjoint from $\bigcup_{j=1}^m C_j$.}

\medskip{\bf Proof.}
We may assume 
 that $D:= \bigcup_{j=1}^m C_j  \neq \emptyset $ and $\|a\|>0$.
In the special case $E=\R$, $\|\cdot\|=|\cdot|$
and $a>0$, we choose $r$ such that $\min D= \min C_r$ if $\min D$ exists,
and   $\inf D= \inf C_r$ otherwise. In  the general case
we apply  the Hahn-Banach theorem (compare e.g.\ Rudin (1991), Theorem 3.3 
and its Corollary) to yield  a linear functional $\ell$ on $E$
with $\ell(a)=\|a\|$ and $|\ell(x)| \le \|x\|$ for every $x\in E$,
so that the special case applied to $\ell(a) > 0$ and
$\ell(C_1),\ldots,\ell(C_m) \subset \R$  yields the  claim.
 \Halmos

\subsection{Lemma}\label{Lemma on Fnm}
{\sl 
For $n,m\in \N_0$, we have  $F_n(m) = \sum_{i=r}^s {n\choose i}$ with 
$r=r_{n,m}:=\lfloor (n-m+1)/2\rfloor$ and $s = r+m-1$,
and also with $\lceil (n-m+1)/2\rceil$ in place of $r_{n,m}$.
Further, 
\la        \label{Fnk recursion}
  F_{n}(m)  &=& F_{n-1}(m-1) + F_{n-1}(m+1) \qquad (n,m\in \N)
\al
and $n \mapsto   2^{-n} F_{n}(m)$ is for every $m\in\N_0$ 
a decreasing function.}\\

{\bf Proof.} 
The claim up to \eqref{Fnk recursion} follows easily from 
the symmetry, monotonicity and recursion properties of
the binomial coefficients. The last claim follows,
since the right hand side of  \eqref{Fnk recursion}
is $\le$ $2F_{n-1}(m)$.
\Halmos\\

We write $\sharp A$ for the cardinality of a set $A$.

\subsection{Theorem (essentially Kleitman's (1970)
Theorem I)}\label{Kleitman theorem}
{\sl Let $n,m\in\N$, $a_1,\ldots,a_n \in E$, 
and $C_1,\ldots,C_m\subset E$ with
\la             \label{Kleitman assumption}
 x,y \in C_j &\Rightarrow& \|x-y\| < \min_{i=1}^n \| a_i \|
\al 
for each $j\in\set{1,\ldots,m}$. Then 
\la                                       \label{Kleitman inequality}
  \sharp \Big\{ I\subset \set{1,\ldots,n}\,:\, \sum_{i\in I} a_i 
 \in \bigcup_{j=1}^m C_j\Big\} &\le & F_{n}(m) 
\al
with equality for $E=\R$, $\|\cdot\|=|\cdot|$, $a_1=\ldots=a_n=1$, 
and $C_j=\set{\lfloor(n-m+1)/2\rfloor + j-1}$.}

\medskip
{\bf Proof:} We consider more generally $n,m\in \N_0$ 
and let $G_n(m)$ denote the supremum of the left hand side of 
\eqref{Kleitman inequality} subject to the stated assumptions
on $a_1,\ldots,a_n$ and $C_1,\ldots,C_m$.
Then 
\la
 G_n(0) &=& 0 \,=\, F_n(0) \qquad (n\in\N_0)   \label{Gnk Anf 1}\\ 
 G_0(m) &= & 1 \,=\, F_0(m) \qquad (m\in\N) \label{Gnk Anf 2}
\al
Let $n,m\in\N$. Given $a_1,\ldots,a_n\in E $ and $C_1,\ldots,C_m \subset E$
with \eqref{Kleitman assumption}, let $a:= a_n$ and choose $r$ according to 
Lemma \ref{Kleitman-Lemma}. Then
the left hand side of 
\eqref{Kleitman inequality} is $\sharp A$ with  
\[
 A &:=& \Big\{\epsilon\in\set{0,1}^n\,:\,
    \sum_{i=1}^n \epsilon_i a_i \in \bigcup_{j=1}^m C_j \Big\} 
  \,\,=\,\,
A_1 \times \set{0} \,\,\cup\,\, A_2\times \set{1} 
         \,\,\cup\,\, A_3\times \set{1} 
\]
where 
\[
 A_1 &:=& \Big\{ \epsilon\in\set{0,1}^{n-1}\,:\, 
    \sum_{i=1}^{n-1}\epsilon_i a_i \in \bigcup_{j=1}^m C_j\Big\}\\
 A_2 &:=& \Big\{ \epsilon\in\set{0,1}^{n-1}\,:\, 
    \sum_{i=1}^{n-1}\epsilon_i a_i \in  C_r -a_n\Big\}\\
 A_3 &:=& \Big\{ \epsilon\in\set{0,1}^{n-1}\,:\, 
    \sum_{i=1}^{n-1}\epsilon_i a_i \in \bigcup_{j\neq r} (C_j- a_n)\Big\}
\]
with $A_1\cap A_2=\emptyset$ and thus  
\[
  \sharp A 
  &\le& \sharp A_1 + \sharp A_2 + \sharp A_3 
  \,\, =\,\,\sharp (A_1 \cup A_2) + \sharp A_3 
  \,\,\le\,\,  G_{n-1}(m+1) + G_{n-1}(m-1)
\]  
Hence we have
\la \label{Kleitman Gnk recursion}
 G_n(m) &\le & G_{n-1}(m-1) + G_{n-1}(m+1) \qquad (n,m\in \N)
\al   
Now   \eqref{Fnk recursion}, \eqref{Gnk Anf 1}, \eqref{Gnk Anf 2}
and \eqref{Kleitman Gnk recursion} together imply
$G_n(k)\le F_n(m)$ for all $n,m\in\N_0$, as was to be shown.
The claim about equality is obvious.
\Halmos\\

We call a random vector $X$  symmetric if it has the same 
law as $-X$. We recall the definitions \eqref{Def BerC} and 
\eqref{Def Fnk, sums of largest binomial coefficients} and put
\[
  \STPC_p &:=& \bigconv_{i=1}^n   \left((1-p_i)\delta_0 +\frac{p_i}2
 (\delta_{-1}+\delta_1)\right) \qquad\qquad (p\in[0,1]^n)
\]

\subsection{Theorem (Kanter's (1976) Lemma 4.2 generalized)}
\label{Theorem (Kanter's (1976) Lemma 4.2 generalized)}
{\sl 
 Let $h\in\,]0,\infty[$, $n,m\in \N$, and $p\in[0,1]^n$.
 Then the supremum of 
 \[
  \P\big(\sum_{i=1}^n X_i \in \bigcup_{j=1}^m C_j\big)
 \]
 taken over all 
 measurable $\R$-vector spaces $E$, 
 measurable seminorms $\|\cdot\|$ on $E$,       
 measurable sets $C_1,\ldots,C_m\subset E$ with 
 \[
   x,y \in C_j &\Rightarrow & \|x-y\| <2h 
 \]
 for every $j\in\set{1,\ldots,m}$,
 and all independent and symmetric $E$-valued random vectors
 $X_i$ with 
 \[ 
  \P(\|X_i\| <h) &\le& 1-p_i \qquad (i=1,\ldots,n)
 \]
 is attained for $E=\R$, $\|\cdot\|= |\cdot |$,
 $C_j=\set{0,h} +(2j-m-1)h$, 
 and the $X_i$  symmetric $\R$-valued with $\P(X_i=0)=1-p_i=1-\P(|X_i|=h)$.
 The value of the supremum is   
 \la                   \label{Kanter Lemma 4.2 sup} 
  \STPC_p ([-m+1,m]) 
  &=& \sum_{k=0}^n 2^{-k}F_k(m) \BerC_p (\set{k})
 \al
}\\

{\bf Remark.} Analytically convenient and sharp upper bounds for the 
quantity in  \eqref{Kanter Lemma 4.2 sup} in the special 
case $m=1$ are provided by   Kanter (1976, Lemma  4.3)
and by Mattner \& Roos (2006). It is an open problem to prove
analogous bounds for $m\ge 2$.

\medskip
{\bf Proof.} We may assume $h=1$.
Let  $n$ etc.\ up to the $X_i$ be as stated and let us put
$\pi_i := 1-\P(\|X_i\|<1)$. We may assume that $X_i =(1-B_i)Y_i+B_iR_iZ_i$
with $B_1,\ldots,B_n,Y_1,\ldots,Y_n,$\\
$Z_1,\ldots,Z_n,R_1,\ldots,R_n$ independent,
$B_i \sim \Ber_{p_i}$, 
$Y_i \sim \P(X_i\in \cdot \given \|X_i\| <1)$ $:=$ the conditional
distribution of $X_i$ given $\|X_i\| <1$,
$Z_i \sim \P(X_i\in \cdot \given \|X_i\| \ge 1)$, and
$\P(R_i=-1)=\P(R_i=1)=1/2$. 
Then, with $B:=(B_1,\ldots,B_n)$, with $Q$ denoting the law of
$(Y,Z):= (Y_1,\ldots,Y_n,Z_1,\ldots,Z_n)$, 
and with $|b|:= \sum_{i=1}^n b_i$, we have  
\la
 &&\P\big(\sum_{i=1}^n X_i \in \bigcup_{j=1}^m C_j\big)\nonumber \\
  &=& \sum_{b\in \set{0,1}^n} \P(B=b)\int
   \underset{\phantom{2^{-|b|}F^{}_{|b|}(m)}\le 2^{-|b|}F^{}_{|b|}(m)}
   {\underbrace{\P\Big(\sum_{i=1}^n b_i\frac{R_i+1}2z_i 
     \in \bigcup_{j=1}^m \frac 12 \big(C_j 
    +\sum_{i=1}^n(b_iz_i-(1-b_i)y_i)\big)\Big)}}
 \d Q(y,z)   \label{KL4.2,Step1}\\
 &\le& \text{R.H.S.\eqref{Kanter Lemma 4.2 sup} with $\pi$ instead of $p$}
   \label{KL4.2,Step2}\\
 &\le& \text{R.H.S.\eqref{Kanter Lemma 4.2 sup}} \label{KL4.2,Step3}
\al
Here the inequality in \eqref{KL4.2,Step1}, and hence \eqref{KL4.2,Step2},
follows from  Theorem \ref{Kleitman theorem}, 
with those $z_i$ with $b_i=1$ playing the role of the $a_i$,
and with $\frac 12 \big(C_j     +\sum_{i=1}^n(b_iz_i-(1-b_i)y_i)\big)$
in place of $C_j$. Inequality \eqref{KL4.2,Step3}
is true since $\N_0\ni k\mapsto 2^{-k}F_k(m)$ is decreasing by Lemma 
\ref{Lemma on Fnm},  and 
$[0,1]^n\ni p\mapsto \BerC_{p}$ is increasing
with respect to the coordinatewise order on $[0,1]^n$ and the 
usual stochastic order. 
In the special case $E=\R$ etc.\ as stated, 
we have $Y_i\sim \delta_0$ and 
may replace the distribution of $Z_i$ by $\delta_1$
in deriving \eqref{KL4.2,Step1}, and hence get equality everywhere. 
\Halmos  

\subsection{Corollary} \label{Corollary to Kanter Lemma 4.2, Version 2}
{\sl Let $0<h\le H <\infty$ with 
$m:=\lceil H/h\rceil < H/h + 1/2$,
$n\in\N$, and $p\in [0,1]^n$. Then the supremum of 
\la         \label{Prob in question}
  \P\big(\sum_{i=1}^n X_i \in \,]-H,H]+a\big)
\al
taken over all independent and symmetric $\R$-valued random 
 variables $X_i$ with 
 \la  \label{Assumption Cor 2.5}
  \P(|X_i|<h) &\le& 1-p_i \qquad(i=1,\ldots,n) 
 \al
and all $a\in\R$, is attained for $\P(X_i=0)=1-p_i =1-\P(|X_i|=h)$ 
and $a= mh -H  $.
The value of the supremum is given in \eqref{Kanter Lemma 4.2 sup}.}

\medskip{\bf Proof.} Given $h,H,m,n,p,X_i$ and $a$ as above, we have
\[
  \eqref{Prob in question} 
  &\le & \P\big(\sum_{i=1}^n X_i \in\, ]-mh,mh]+b\big) \\
  &=& \P\big(\sum_{i=1}^n X_i \in \bigcup_{j=1}^m C_j \big) \\
  &\le & \text{R.H.S.\eqref{Kanter Lemma 4.2 sup}}
\]  
with $b:=a$ and $C_j := \,]-h,h] + (2j-m-1)h+b$, using Theorem 
\ref{Theorem (Kanter's (1976) Lemma 4.2 generalized)}
with $E=\R$ and $\|\cdot\| = |\cdot|$. On the other 
hand, if $\P(X_i=0)=1-p_i =1-\P(|X_i|=h)$ 
and $a= mh -H  $, and if we let $b:= 0$ instead of $b:=a$, 
then we can replace the two inequalities in the above calculation by
equalities, as the assumption $m<H/h+1/2$ yields $-mh \le -H+a < -(m-1)h$. 
\Halmos

\subsection{Corollary} \label{Corollary improving Nagaev}
{\sl Let $S=\sum_{i=1}^n X_i$ with independent and symmetric $\R$-valued
random variables $X_i$ and let $h\in\,]0,\infty[$. Then 
\eqref{Improvement of Nagaev Thm 1}  holds with $p_i :=  \P(|X_i| \ge h)$
for $i=1,\ldots,n$.}

\medskip{\bf Proof.} For $t>0$, we apply Corollary  
\ref{Corollary to Kanter Lemma 4.2, Version 2} with $a=0$ and 
$H=mh$ with 
$m:= \lfloor  t/h \rfloor +1$ to get
\[
 \P(|S| \le t) &\le& \P(S \in \,]-mh,mh]) \,\,\le \,\, 
 \text{R.H.S.\eqref{Kanter Lemma 4.2 sup}}
\] 
Inequality  \eqref{Improvement of Nagaev Thm 1} follows by taking complements,
since  $F_k(m)=2^k$ for $k\le m -1$. \Halmos\\

\section{Comparison theorems}\label{Comparison theorems} \mbox{} \\

For $\R$-valued random variables $U$ and $V$, 
we write $U \stge V$ if $U$ is stochastically
larger than $V$, that is, if $\P(U \ge t) \ge \P(V\ge t)$ for every 
$t\in\R$.  A specialization of
Corollary \ref{Corollary to Kanter Lemma 4.2, Version 2} 
can be viewed as one of three results yielding at least almost a stochastic
ordering $|S| \stge |T|$ for  sums $S,T$  of independent symmetric 
random variables assuming a corresponding ordering of their terms, the other
two results being theorems of Pruss (1997) and Birnbaum (1948).
It therefore appears natural to summarize these results here,  
and to use  this opportunity to extend Birnbaum's theorem to the lattice 
case.

Let us agree on the following unimodality definitions for laws 
$P$ on $\R$. We call $P$  
{\em unimodal on $\R$}, if $P$ is unimodal in the usual sense that,
for some $x_0\in \R$, the distribution
function of $P$ is convex on $]-\infty,x_0[$ and concave 
on $]x_0,\infty[$.  For $a\in \R$ and $h\in\,]0,\infty[$, we call
$P$  {\em unimodal on $h \Z+a$}, if
$P(h\Z+a)=1$ and if there is a $k_0\in\Z$ such that 
$k\mapsto P(\set{hk+a})$ is increasing on $\set{k\in\Z: k\le k_0}$
and decreasing on $\set{k\in\Z: k\ge k_0}$. 
For $h\in [0,\infty[$, we call $P$ {\em unimodal with span $h$},
if either $h=0$ and $P$ is unimodal on $\R$, or $h>0$ and
$P$ is unimodal on $h\Z+a$ for some $a\in\R$. 
As usual, we attribute any property just defined to a random 
variable $X$ if its distribution enjoys it.

\subsection{Theorem} \label{Comp thm}
{\sl Let $n\in\N$ and let $X_1,\ldots,X_n$ 
as well as $Y_1,\ldots,Y_n$ be independent and symmetrically 
distributed $\R$-valued random variables with 
sums $S=\sum_{i=1}^n X_i$ and $T=\sum_{i=1}^n Y_i$ and with 
\la  \label{|Xi| stgr |Yi|}
   |X_i| &\stge&  |Y_i|\qquad\qquad (i=1,\ldots,n)
\al
\subsubsection{(Pruss (1997))} \label{Pruss' Theorem} 
Then
\[
  \P( |S| \ge t) &\ge&  \frac 12 \P(|T|\ge t)\qquad\qquad (t>0)   
\]
\subsubsection{} \label{Cor to Cor 2.5}
If $h\in\,]0,\infty[$ and $\P(Y_i\in\set{-h,0,h})=1$
for $i=1,\ldots,n$, then 
\la                \label{|S| almost stge |T|}
 \P( |S| > mh )+\frac 12 \P(|S| =mh) &\ge &
   \P( |T| > mh )+\frac 12 \P(|T| =mh)  \qquad (m\in\N)
\al
\subsubsection{(Birnbaum (1948) generalized)} \label{Birnbaum gen}
Let $h\in [0,\infty[$ and $X_1,\ldots,X_n,Y_1,\ldots,Y_n$ be unimodal
with span $h$. In case of $h>0$ assume further for each
$i\in\set{1,\ldots,n}$ that $X_i,Y_i$ are both 
$h\Z$-valued or both $h(\Z +\frac 12)$-valued. 
Then $|S| \stge |T|$.}\\

See Berger (1997, Theorem 1.1) for a further related comparison theorem.\\ 

{\bf Example.} Let 
$n=2$, $X_1,X_2,Y_1 \sim \frac 12 (\delta_{-1}+\delta_1)$, and 
$Y_2=0$. Then $|X_i|\stge |Y_i|$ for $i=1,2$. Since    
$ \P( |S| \ge 1) \,=\, \frac 12$ 
and $\P(|T\ge 1)=1$, 
it follows that the constant $\frac 12$ in Pruss' theorem is best possible.
As each of the four random variables is unimodal with span $2$, 
it also follows that the second sentence in part
\ref{Birnbaum gen} can not be omitted.
Further, in this example, $\P(Y_i\in\set{-1,0,1})=1$ for $i=1,2$ but 
$\P(S>0)+\frac 12 \P(S=0) = \frac 34 \not \ge 1
=\P(T>0)+\frac 12 \P(T=0)$,
showing that in  
\eqref{|S| almost stge |T|} we may not replace $\N$ by $\N_0$.

\medskip{\bf Proof.} {\bf\ref{Pruss' Theorem}} 
See Pruss (1997).

{\bf \ref{Cor to Cor 2.5}} Here \eqref{|Xi| stgr |Yi|} is equivalent 
to \eqref{Assumption Cor 2.5} with $p_i=\P(|Y_i|=h)$, so that 
Corollary \ref{Corollary to Kanter Lemma 4.2, Version 2} with
$H=mh$ and $a=0$ 
yields \eqref{|S| almost stge |T|}.
 
{\bf \ref{Birnbaum gen}} Induction based on Lemmas
\ref{Birmbaum for two} and \ref{Wintner} given below.
In the step from $n-1$ to $n$, we may assume $X_1,\ldots,X_n,Y_1,\ldots,Y_n$
to be independent, and conclude that 
\[
  \Big|\sum_{i=1}^n X_i \Big| &\stge& \Big|\sum_{i=1}^{n-1} X_i \,+\, Y_n\Big|
   \,\,\stge\,\,  \Big|\sum_{i=1}^n Y_i \Big|
\]
by applying  Lemma \ref{Birmbaum for two} first to 
$U_1:= \sum_{i=1}^{n-1} X_i$, $V_1 := X_n$, $W_1 := Y_n$
and then to $U_2:= Y_n$, $V_2:= \sum_{i=1}^{n-1} X_i$, 
$W_2:= \sum_{i=1}^{n-1} Y_i$, observing that by Lemma \ref{Wintner}
the sum $U_1$ is symmetric and unimodal with span $h$,
and that in case of $h>0$ the sums $V_2,W_2$ are both $h\Z$-valued
or both $h(\Z+\frac 12)$-valued.
\Halmos

\subsection{Lemma}\label{Birmbaum for two}
{\sl Let $U,V,W$ be symmetrically distributed 
$\R$-valued random variables with $U,V$ independent, $U,W$ independent,
and $|V| \stge |W|$. 
Let $h\in[0,\infty[$ with $U$ unimodal with span $h$.
In case of $h>0$ let further  $V,W$ be  both 
$h\Z$-valued or both $h(\Z +\frac 12)$-valued.  
Then $|U+V| \stge |U+W|$.}\\

{\bf Proof.} We may assume that $h\in \set{0,1}$.
In case of $h=0$ we put $A:=B:= [0,\infty[$, while for $h=1$
we  let $A,B\in\set{\N_0,\N_0+\frac 12}$
with $\P(|U|\in A)= \P(|V|\in B)=1$.  Then for 
$t\in A+B := \set{a+b\,:\,a\in A,\, b\in B}$
and denoting by $P^{}_U$ etc.\ the laws of the random variables 
occuring as subscripts, we  have
\[
 \P(|U+V| \le t) &=& \int_{B} P^{}_U([v-t,v+t]) \d P^{}_{|V|}(v) \\
  &\le& \int_{B} P^{}_U([v-t,v+t]) \d P^{}_{|W|}(v) \\
  &=&  \P(|U+W| \le t)
\]
since in each case the function $B\ni v \mapsto  P^{}_U([v-t,v+t]) $
is decreasing. As $\P(|U+V| \in A+B) = 1$, this  proves 
$|U+V| \stge |U+W|$.\Halmos

\subsection{Lemma (Wintner)}\label{Wintner}
{\sl Let $X$ and $Y$ be independent $\R$-valued random
variables and let $h\in[0,\infty[$. If $X$ and $Y$ are symmetric
and unimodal with span $h$, then so is $X+Y$.}\\

{\bf Proof.} Obvious by writing the laws of $X$ and $Y$ 
as mixtures of uniform distributions on  symmetric intervals
in $\R$ or $h\Z$ or $h(\Z+\frac 12)$.
See  Dharmadhikari \& Joag-Dev (1988, pp. 13 and 109)
for the cases where  $h=0$ or $X$ and $Y$ 
are both symmetric unimodal on $h\Z$. The remaining 
three cases are analogous. \Halmos

\newpage
\section{Historical notes}\label{Historical notes}
 \mbox{} \\

Theorem \ref{Kleitman theorem} in the Hilbert space case,
and  assuming the sets $C_j$ to be slightly smaller than necessary, 
was proved by Kleitman (1970), generalizing several earlier results
and in particular the one-dimensional case due to 
Erd\H{o}s (1945, Theorems 1 and 3).
Jones (1978, page 4, footnote 7) observed that Kleitman's result and 
proof extends to general (semi-)normed spaces. 
Meanwhile,  Kanter (1976, Lemma 4.1) proved a weaker result,
assuming in particular symmetry of the sets $C_j$.
The present proof of Theorem \ref{Kleitman theorem} is just
a slightly refined rewrite  of Kleitman's proof and Jones' footnote.

Kanter (1976) essentially stated and proved
Theorem \ref{Theorem (Kanter's (1976) Lemma 4.2 generalized)}
for $m=1$ and $C_1$ symmetric. Le Cam (1986, pp.\ 408-409)
adopted Kanter's approach. 

Theorem \ref{Comp thm}\ref{Birnbaum gen} in the case 
of $h=0$ and without atoms at zero is due to 
Birnbaum (1948).
Bickel \& Lehmann (1976) and Shaked \& Shantikumar (1994, 
page 78) allowed  atoms at zero
in their statements, but apparently not in their proofs.
Sherman (1955) extended Birnbaum's result to the absolutely continuous
multivariate case. Dharmadhikari \& Joag-Dev (1988, p.\ 164)
gave an elegant development of Sherman's theorem, dispensing 
with unnecessary continuity assumptions. 
They also essentially stated without proof Theorem 
\ref{Comp thm}\ref{Birnbaum gen} for $h>0$ in the case where 
all random variables are $h\Z$-valued.

\bigskip
{\Large\bf References } 

\medskip
{\footnotesize
\begin{description}
\item[\sc Berger, E. (1997).] Comparing sums of independent bounded 
random variables and sums of Bernoulli random variables.
\SPL {\bf 34}, 251-258.

\item[\sc Bickel, P.J. \& Lehmann, E.L. (1976).] Descriptive statistics 
for nonparametric  models. III. Dispersion. 
{\it Ann.\ Statist.} {\bf 4}, 1139-1158.

\item[\sc Birnbaum, Z.W. (1948).] 
On random variables with comparable peakedness. 
{\it Ann.\ Math.\ Statist.} {\bf 19}, 76-81.

\item[\sc Dharmadhikari, S. \& Joag-Dev, K. (1988).]
{\it Unimodality, Convexity, and Applications.} Academic Press, San Diego.

\item[\sc Erd\H{o}s, P. (1945).] On a lemma of Littlewood and Offord.
\BAMS {\bf 51}, 898-902.

\item[\sc Jones, L. (1978).] 
On the distribution of sums of vectors.
{\it SIAM J.\ Appl.\ Math.} {\bf 34}, 1-6.

\item[\sc Le Cam (1986).] {\it Asymptotic Methods in Statistical
Decision Theory.} Springer-Verlag, New York.

\item[\sc Kanter, M. (1976).] Probability inequalities for convex sets
and multidimensional concentration functions.
{\it J. Multivariate Anal.} {\bf 6}, 222-236.  

\item[\sc Kleitman, D. (1970).] On a lemma of Littlewood and Offord
on the distributions of linear combinations of vectors.
{\it Advances in Math.} {\bf 5}, 155-157.

\item[\sc Mattner, L. \& Roos, B. (2006).] 
A shorter proof of Kanter's Bessel function concentration bound.
Preprint. Available at arXiv.math.PR/0603522

\item[\sc Nagaev, S.V. (2001).] Lower bounds for probabilities of large
  deviations of sums of independent random variables.
 \TPA {\bf 46}, 728-735. 

\item[\sc Pruss, A.R. (1997).] Comparisons between tail probabilities of sums
of independent symmetric random variables. {\it Ann.\ Inst.\ Henri Poincar\'e}
{\bf 33}, 651-671.
\item[\sc Rudin, W. (1991).] {\it Functional Analysis.}
2nd ed. McGraw-Hill, N.Y.

\item[\sc Shaked, M. \& Shantikumar, J.G. (1994).]
{\it Stochastic Orders and their Applications.}
Academic Press, San Diego.

\item[\sc Sherman, S. (1955).]
A theorem on convex sets with applications.
\AnnMS {\bf 26}, 763-767.

\end{description}
}
{\footnotesize
\noindent
{\sc Universit\"at zu L\"ubeck\\ 
     Institut f\"ur Mathematik\\
     Wallstr.\ 40\\
     D-23560 L\"ubeck\\
     Germany\\
     Email: \verb§mattner@math.uni-luebeck.de§ 
}  
}
\end{document}